# Time Series and Related Topics

**In Memory of Ching-Zong Wei**

**Hwai-Chung Ho, Ching-Kang Ing, Tze Leung Lai, Editors**

Institute of Mathematical Statistics
*Lecture Notes–Monograph Series*

Series Editor:
R. A. Vitale

Printed in the United States of America

# Contents

iii

RELATED TOPICS

# Contributors to this volume

Aston, J. A. D. *Academia Sinica*

Breidt, F. J. *Colorado State University*

Cantor, J. L. *Science Application International Corporation*
Chan, H. P. *National University of Singapore*
Chan, N. H. *The Chinese University of Hong Kong*
Chen, L. H. Y. *National University of Singapore*
Chiang, T.-S. *Academia Sinica*

Davis, R. A. *Colorado State University*

Findley, D. F. *U.S. Census Bureau*
Fuh, C.-D. *Academia Sinica*

Giurcăneanu, C. D. *Tampere University of Technology*

Ho, H.-C. *Academia Sinica*
Hsiao, C. *University of Southern California*
Hsu, N.-J. *National Tsing-Hua University*
Hu, I. *Hong Kong University of Science and Technology*

Ing, C.-K. *Academia Sinica*

Lai, T. L. *Stanford University*
Li, W. K. *The University of Hong Kong*
Lin, J.-L. *Academia Sinica*

Ng, C. T. *The Chinese University of Hong Kong*
Nielsen, B. *University of Oxford*

Pötscher, B. M. *University of Vienna*

Rissanen, J. *Technical University of Tampere and Helsinki, and Helsinki Institute for Information Technology*
Robinson, P. M. *London School of Economics*
Rosenblatt, M. *University of California at San Diego*

Sheu, S.-J. *Academia Sinica*
Shiu, S.-Y. *University of Utah*
Sin, C.-Y. *Xiamen University*

Tsay, R. S. *University of Chicago*

Wei, C.-Z. *Academia Sinica*
Wong, S. P.-S. *The Chinese Universty of Hong Kong*
Wong, T. S. T. *The University of Hong Kong*

Xia, A. *University of Melbourne*

Yao, Y.-C. *Academia Sinica*

# Preface

A major research area of Ching-Zong Wei (1949–2004) was time series models and their applications in econometrics and engineering, to which he made many important contributions. A conference on time series and related topics in memory of him was held on December 12–14, 2005, at Academia Sinica in Taipei, where he was Director of the Institute of Statistical Science from 1993 to 1999. Of the forty-two speakers at the conference, twenty contributed to this volume. These papers are listed under the following three headings.

## 1. Estimation and prediction in time series models

Breidt, Davis, Hsu and Rosenblatt consider estimation of the unknown moving average parameter $\theta$ in an MA(1) model when $\theta = 1$, and derive the limiting pile-up probabilities $P(\hat{\theta} = 1)$ and $1/n$-asymptotics for the Laplace likelihood estimator $\hat{\theta}$. Cantor and Findley introduce a recursive estimator for $\theta$ in a possibly misspecified MA(1) model and obtain convergence results by approximating the recursive algorithm for the estimator by a Robbins–Monro-type stochastic approximation scheme. Giurcăneanu and Rissanen consider estimation of the order of AR and ARMA models by stochastic complexity, which is the negative logarithm of a normalized maximum likelihood universal density function. Nielsen investigates estimation of the order in general vector autoregressive models and shows that likelihood-based information criteria, and likelihood ratio tests and residual-based tests can be used, regardless of whether the characteristic roots are inside, or on, or outside the unit disk, and also in the presence of deterministic terms. Instead of model selection, Pötscher considers model averaging in linear regression models, and derives the finite-sample and asymptotic distributions of model averaging estimators. Robinson derives the asymptotic properties of conditional-sum-of-squares estimates in parametric models of stationary time series with long memory. Ing and Sin consider the final prediction error and the accumulated prediction error of the adaptive least squares predictor in stochastic regression models with nonstationary regressors. The paper by Lin and Wei, which was in preparation when Ching-Zong was still healthy, investigates the adaptive least squares predictor in unit-root nonstationary processes.

## 2. Time series modeling in finance, macroeconomics and other applications

Aston considers criteria for deciding when and where heavy-tailed models should be used for macroeconomic time series, especially those in which outliers are present. Hsiao reviews nonstationary time series analysis from the perspective of the Cowles Commission structural equation approach, and shows that the same rank condition for identification holds for both stationary and nonstationary time series, that certain instrumental variables are needed for consistent parameter estimation, and that classical instrumental-variable estimators have to be modified for valid inference in the presence of unit roots. Chan and Ng investigate option pricing when the volatility of the underlying asset follows a fractional version of the CEV (constant elasticity of variance) model. Ho considers linear process models, with a latent long-memory volatility component, for asset returns and provides asymptotically normal estimates, with a slower convergence rate than $1/\sqrt{n}$, of the Sharpe ratios in these investment models. Tsay reviews some commonly used models for the time-varying multivariate volatility of $k$ ($\geq 2$) assets and proposes a simple parsimonious approach that satisfies positive definite constraints on the time-varying correlation matrix. Lai and Wong propose a new approach to time series modeling that combines subject-matter knowledge of the system dynamics with statistical techniques in time series analysis and regression, and apply this approach to American option pricing and the Canadian lynx data.

## 3. Related topics

Besides time series analysis, Ching-Zong also made important contributions to the multi-armed bandit problem, estimation in branching processes with immigration, stochastic approximation, adaptive control and limit theorems in probability, and had an active interest in the closely related areas of experimental design, stochastic control and estimation in non-regular and non-ergodic models. The paper by Chan, Fu and Hu uses the multi-armed bandit problem with precedence relations to analyze a multi-phase management problem and thereby establishes the asymptotic optimality of certain strategies. Yao develops an approximation to Gittins index in the discounted multi-armed bandit problem by using a continuity correction in an associated optional stopping problem. Chen and Xia describe Stein's method for Poisson approximation and for Poisson process approximation from the points of view of immigration-death processes and Palm distributions. Cheng, Wu and Huwang propose a new approach, which is based on a response surface model, to the analysis of experiments that use the technique of sliding levels to treat related factors, and demonstrate the superiority of this approach over previous methods in the literature. Chiang, Sheu and Shiu formulate the valuation problem of a financial derivative in markets with transaction costs as a stochastic control problem and consider optimization of expected utility by using the price systems for these markets. Wong and Li propose to use the maximum product of spacings (MPS) method for parameter estimation in the GEV (generalized extreme value) family and the generalized Pareto family of distributions, and show that the MPS estimates are asymptotically efficient and can outperform the maximum likelihood estimates.

We thank the Institute of Statistical Science of Academia Sinica for providing financial support for the conference. Special thanks also go to the referees who reviewed the manuscripts. A biographical sketch of Ching-Zong and a bibliography of his publications appear after this Preface.

<div style="text-align: right">

Hwai-Chung Ho
Ching-Kang Ing
Tze Leung Lai

</div>

# Biographical sketch

Ching-Zong Wei was born in 1949 in south Taiwan. He studied mathematics at National Tsing-Hua University, Taiwan, where he earned a BS degree in 1971 and an MS degree in 1973. He went to the United States in 1976 to pursue advanced studies in statistics at Columbia University, where he earned a PhD degree in 1980. He then joined the Department of Mathematics at the University of Maryland, College Park, as an Assistant Professor in 1980, and was promoted to Associate Professor in 1984 and Full Professor in 1988. In 1990 he returned to Taiwan, his beloved homeland, to join the Institute of Statistical Science at Academia Sinica, where he stayed as Research Fellow for the rest of his life, serving between 1993 and 1999 as Director of the Institute. He also held a joint appointment with the Department of Mathematics at National Taiwan University.

In addition to his research and administrative work at Academia Sinica, Ching-Zong also made important contributions to statistical education in Taiwan. To promote statistical thinking among the general public, he published in local newspapers and magazines articles on various topics of general interest such as lottery games and the Bible code. These articles, written in Chinese, introduced basic statistical and probabilistic concepts in a heuristic and reader-friendly manner via entertaining stories, without formal statistical jargon.

Ching-Zong made fundamental contributions to stochastic regression, adaptive control, nonstationary time series, model selection and sequential design. In particular, his pioneering works on (i) strong consistency of least squares estimates in stochastic regression models, (ii) asymptotic behavior of least squares estimates in unstable autoregressive models, and (iii) predictive least squares principles in model selection, have been influential in control engineering, econometrics and time series. A more detailed description of his work appears in the Bibliography. He was elected Fellow of the Institute of Mathematical Statistics in 1989, and served as an Associate Editor of the *Annals of Statistics* (1987–1993) and *Statistic Sinica* (1991–1999). In 1999, when Ching-Zong was at the prime of his career, he was diagnosed with brain tumors. He recovered well after the first surgery and remained active in research and education. In 2002, he underwent a second surgery after recurrence of the tumors, which caused deterioration of his vision. He continued his work and courageous fight with brain tumors and passed away on November 18, 2004, after an unsuccessful third surgery. He was survived by his wife of close to 30 years, Mei, and a daughter. In recognition of his path-breaking contributions, Vol. 16 of *Statistica Sinica* contains a special memorial section dedicated to him.

# Bibliography

Before listing Ching-Zong's publications, we give a brief introduction of their background and divide them broadly into five groups, in which the papers are referred to by their numbers in the subsequent list.

## A. Least squares estimates in stochastic regression models

Ching-Zong's work in this area began with papers [1], [2] and [3], in which the strong consistency of least squares estimates is established in fixed-design linear regression models. In particular, when the errors are square integrable martingale differences, a necessary and sufficient condition for the strong consistency of least squares estimates is given. However, when the regressors are stochastic, this condition is too weak to ensure consistency. Paper [6] is devoted to resolving this difficulty, and establishes strong consistency and asymptotic normality of least squares estimates in stochastic regression models under mild assumptions on the stochastic regressors and errors. These results can be applied to interval estimation of the regression parameters and to recursive on-line identification and control schemes for linear dynamic systems, as shown in [6]. Papers [7], [12] and [15] extend the results of [6] and establish the asymptotic properties of least squares estimates in more general settings.

## B. Adaptive control and stochastic approximation

Papers [17] and [18] resolve the dilemma between the control objective and the need of information for parameter estimation by occasional use of white-noise probing inputs and by a reparametrization of the model. Asymptotically efficient self-tuning regulators are constructed in [18] by making use of certain basic properties of adaptive predictors involving recursive least squares for the reparametrized model. Paper [16] studies excitation properties of the designs generated by adaptive control schemes. Instead of using least squares, [13] uses stochastic approximation for recursive estimation of the unknown parameters in adaptive control. Paper [20] introduces a multivariate version of adaptive stochastic approximation and demonstrates that it is asymptotically efficient from both the estimation and control points of view, while [28] uses martingale transforms with non-atomic limits to analyze stochastic approximation. Paper [23] introduces irreversibility constraints into the classical multi-armed bandit problem in adaptive control.

## C. Nonstationary time series

For a general autoregressive (AR) process, [9] proves for the first time that the least squares estimate is strongly consistent regardless of whether the roots of the characteristic polynomial lie inside, on, or outside the unit disk. Paper [22] shows that in general unstable AR models, the limiting distribution of the least squares estimate can be characterized as a function of stochastic integrals. The techniques developed in [22] and in the earlier paper [19] for deriving the asymptotic distribution soon became standard tools for analyzing unstable time series and led to many important developments in econometric time series, including recent advances in the analysis of cointegration processes.

### D. Adaptive prediction and model selection

Paper [21] considers sequential prediction problems in stochastic regression models with martingale difference errors, and gives an asymptotic expression for the cumulative sum of squared prediction errors under mild conditions. Paper [27] shows that Rissanen's predictive least squares (PLS) criterion can be decomposed as a sum of two terms; one measures the goodness of fit and the other penalizes the complexity of the selected model. Using this decomposition, sufficient conditions for PLS to be strongly consistent in stochastic regression models are given, and the asymptotic equivalence between PLS and the Bayesian information criterion (BIC) is established. Moreover, a new criterion, FIC, is introduced and shown to share most asymptotic properties with PLS while removing some of the difficulties encountered by PLS in finite-sample situations. In [38], the first complete proof of an analogous property for Akaike's information criterion (AIC) in determining the order of a vector autoregressive model used to fit a weakly stationary time series is given, while in [41], AIC is shown to be asymptotically efficient for same-realization predictions. Closely related papers on model selection and adaptive prediction are [39], [42] and [43].

### E. Probability theory, stochastic processes and other topics

In [4] and [5], sufficient conditions are given for the law of the iterated logarithm to hold for random subsequences, least squares estimates in linear regression models and partial sums of linear processes. Papers [8] and [14] provide sufficient conditions for a general linear process to be a convergence system, while [10] considers martingale difference sequences that satisfy a local Marcinkiewicz-Zygmund condition. Papers [24], [25] and [26] resolve long-standing estimation problems in branching processes with immigration. Paper [35] studies the asymptotic behavior of the residual empirical process in stochastic regression models. In [36], uniform convergence of sample second moments is established for families of time series arrays, whose modeling by multistep prediction or likelihood methods is considered in [40]. Paper [11], [29], [30] and [33] investigate moment inequalities and their statistical applications. Density estimation, mixtures, weak convergence of recursions and sequential analysis are considered in [31], [32], [34] and [37].

### Publications of Ching-Zong Wei

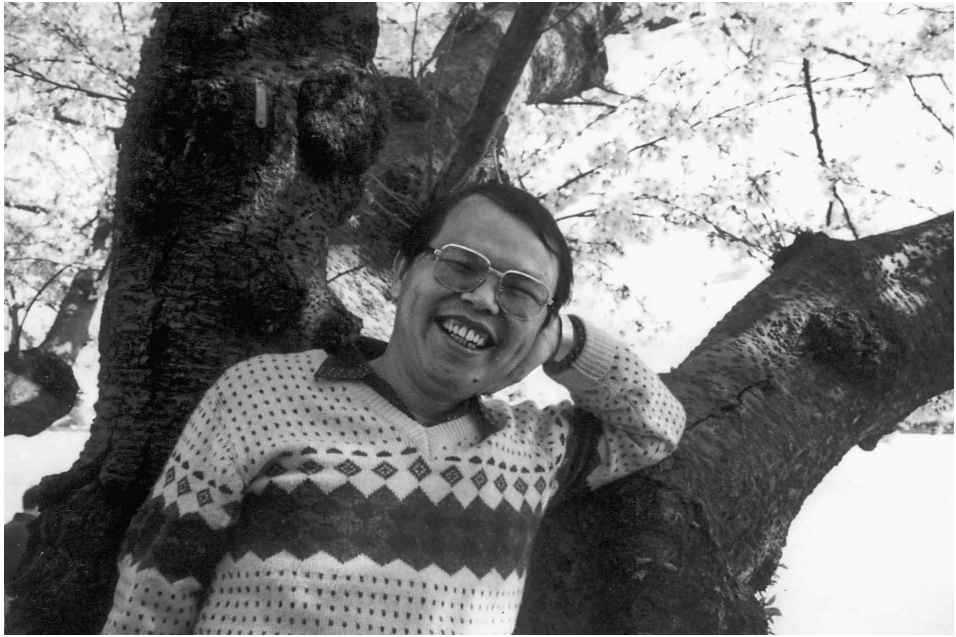

Ching-Zong Wei, Maryland 1985.

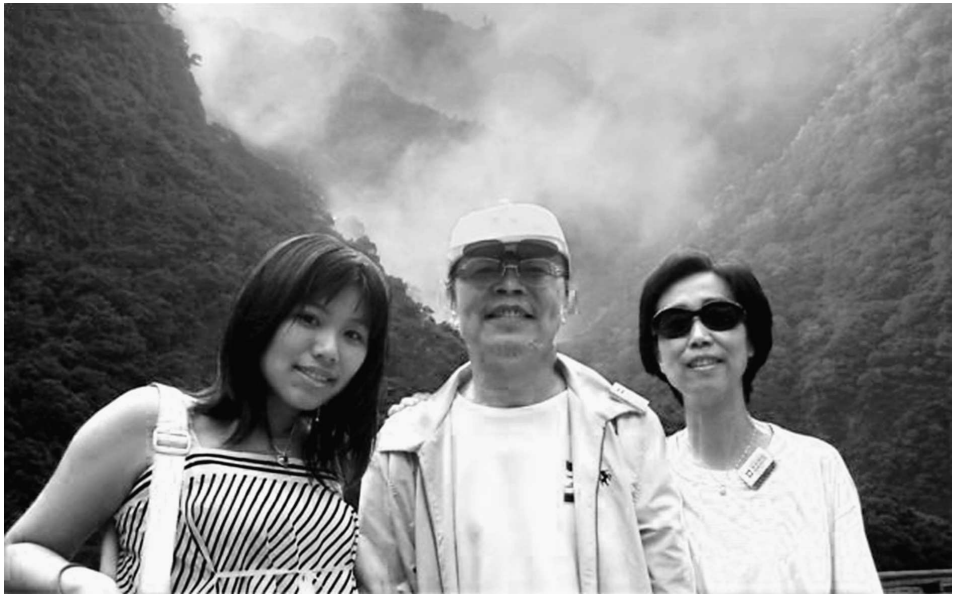

In Hualian, Taiwan, with wife and daughter, 2004.